\documentclass[12pt,letterpaper]{article}
\usepackage[T1]{fontenc}
\usepackage{palatino}
\usepackage{mathpazo}
\usepackage{amssymb, amstext}
\usepackage{amsmath, amsthm, fullpage}
\usepackage{bbm}
\usepackage{graphicx}
\usepackage{calc}
\usepackage{caption}
\newcommand{\R}{\mathbbm{R}}

\DeclareMathOperator{\codim}{codim}

\def\eqbd{\mathop{{:}{=}}}

\newcommand{\norm}[2]{\lVert{#1}\rVert_{#2}}

\newcommand{\abs}[1]{\lvert{#1}\rvert}

\newcommand{\const}{{\rm const}}

\newtheorem{theorem}{Theorem}

\newtheorem*{definition}{Definition}
\newtheorem*{observation*}{Observation}
\newtheorem*{assume*}{Assumption}
\title{Compressive sensing: \\ a paradigm shift in signal processing}

\author{Olga V. Holtz \\ \small University of California-berkeley \\ \small \& Technische Universit\"at Berlin}

\date{}

\begin{document}
\maketitle

\abstract{\small \noindent We survey a new paradigm in signal processing known as "compressive 
sensing". Contrary to old practices of data acquisition and reconstruction based 
on the Shannon-Nyquist sampling principle, the new theory shows that it is 
possible to reconstruct images or signals of scientific interest accurately and  
even exactly from a number of samples which is far smaller than the desired 
resolution of the image/signal, e.g., the number of pixels in the image. This
new technique draws from results in several fields of mathematics, including
algebra, optimization, probability theory, and harmonic analysis. We will 
discuss some of the key mathematical ideas behind compressive sensing, as 
well as its implications to other fields: numerical analysis, information 
theory, theoretical computer science, and engineering.}



\section{Introduction}

{\it Compressive sensing}~\cite{Don06, TsaD06a} is a new concept in signal processing where one seeks to
minimize the number of measurements to be taken from signals while still retaining
the information necessary to approximate them well. The ideas have their origins
in certain abstract results from functional analysis and approximation theory~\cite{Kashin77,Pin85} 
but were recently brought into the forefront by the work of Cand\'{e}s, Romberg and Tao~\cite{CanR06,CanT06,CanRT06} 
and Donoho~\cite{Don06} who constructed concrete algorithms and
showed their promise in application. 

Sparse approximation has been studied for nearly a century, and it has numerous applications.
Temlyakov \cite{Tem02} locates the first example in a 1907 paper of Schmidt \cite{Sch07}. In the 1950s,
statisticians launched an extensive investigation of another sparse approximation problem called
subset selection in regression \cite{Mil02} and recently  least angle regression \cite{EfrHJT04, Tib96}. Later, approximation theorists began a systematic study of $m$-term approximation with respect to orthonormal bases and redundant systems \cite{DeV98, Tem02} and very recently in \cite{CohDD06,CohDD07}.

Over the last decade, the signal processing community spurred by the work of Coifman et al. \cite{CoiM89, CoiW92} and Mallat et al. \cite{MalZ93, DavMZ94, DavMA97} has become interested in sparse representations
for compression and analysis of audio \cite{GriB03}, images \cite{FroVVK04} and video \cite{NguZ03}.
Sparsity criteria also arise in deconvolution \cite{TayBM79}, signal modeling \cite{Ris79}, pre-conditioning
\cite{GroH97}, machine learning \cite{Gir98}, de-noising \cite{CheDS99}, regularization \cite{DauDD04, DauTV07} and error correction \cite{CanT06a, CanT05, Fel02, Fel03a, Fel03b, Fel03c}.
Most sparse approximation problems employ a linear model in which the collection of elementary
signals is both linearly dependent and large. These models are often called redundant or overcomplete.
Recent research suggests that overcomplete models offer a genuine increase in approximation
power \cite{RaoB98, FroV03}. Unfortunately, they also raise a serious challenge. How do we find a good
representation of the input signal among the plethora of possibilities? One method is to select a
parsimonious or sparse representation. The exact rationale for invoking sparsity may range from
engineering to economics to philosophy. At least three justifications are commonly given:
\begin{enumerate}
\item It is sometimes known \emph{a priori} that the input signal can be expressed as a short linear
combination of elementary signals also contaminated with noise.
\item The approximation may have an associated cost that must be controlled. For example, the
computational cost of evaluating the approximation depends on the number of elementary
signals that participate. In compression, the goal is to minimize the number of bits required
to store the approximation.
\item Some researchers cite Occam's Razor, ''Pluralitas non est ponenda sine necessitate'' (causes
must not be multiplied beyond necessity).
\end{enumerate}

Sparse approximation problems are computationally challenging because most reasonable sparsity
measures are not convex. A formal hardness proof for one important class of problems independently
appeared in \cite{Nat95} and \cite{DavMA97}. A vast array of heuristic methods for producing sparse approximations 
have been proposed, but the literature contains few guarantees of their performance. The pertinent numerical 
techniques fall into at least three basic categories:
\begin{enumerate}
\item The convex relaxation approach replaces the nonconvex sparsity measure with a related convex function to obtain a convex programming problem. The convex program can be solved in polynomial time with standard software \cite{BoyV04}, and one expects that it will yield a good sparse approximation. More on that will be said in the sequel.
\item  Greedy methods make a sequence of locally optimal choices in an effort to produce a good
global solution to the approximation problem. This category includes forward selection
procedures (such as matching pursuits), backward selection and others. Although these
approaches sometimes succeed \cite{CouB00, GilGIMS02, GilMS03, GilMStr05, Tro04, TroG05, TroGMS03}, they can also fail spectacularly
\cite{DeVT96, CheDS99}. The monographs of Miller \cite{Mil02} and Temlyakov \cite{Tem02} taste
the many flavors of greedy heuristic.
\item Specialized nonlinear programming software has been developed that attempts
to solve sparse approximation problems directly using, for example, interior point methods
\cite{RaoK99}. These techniques are only guaranteed to discover a locally optimal solution though.
\end{enumerate}
\bigskip
Several problems require solutions to be obtained from underdetermined systems of linear equations, i.e., systems with fewer equations than unknowns. Some example of such problems arise in linear filtering signal processing, and inverse problems. For an underdetermined system of linear equations, if there is any solution, there are infinitely many solutions. In many applications, the ``simplest`` solution is most acceptable. Such a solution is inspired by the minimalist principle of Occam's Razor. For example, if the parameters of a model are being estimated then among all models that explain the data equally well, the one with the minimum number of parameters is most desirable.

\begin{figure}\label{Loganshepp1}
\begin{center}
\includegraphics[angle=0, width=0.6\textwidth]{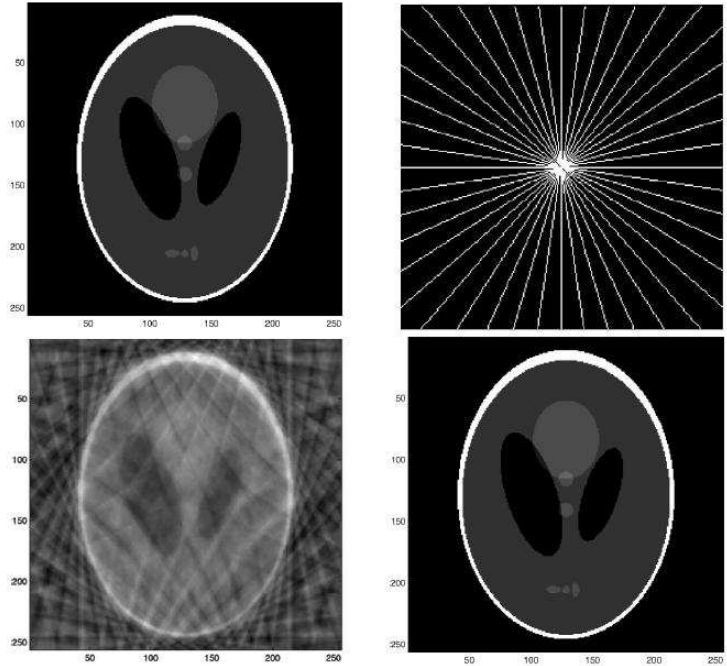}
\end{center}  
\caption{When Fourier coefficients of a testbed medical image known as
the Logan--Shepp phantom (top left) are sampled along $22$ radial lines 
in the frequency domain (top right), a naive, ``minimal energy'' reconstruction 
setting unobserved Fourier coefficients to $0$ is marred by artifacts (bottom left). 
$\ell_1$-reconstruction (bottom right) is exact.}
\end{figure}

The notion that {\it sparse} signals--meaning signals with a small number of nonzero coefficients for a given basis (and no noise) one can (with high probability) be reconstructed exactly via $\ell_1$-minimization is not exactly new. The idea was first expressed in 1986, by Fadil Santosa and William Symes~\cite{SanS86}. But the full extent of the theory, including the robustness of the reconstruction procedure, is only now coming into full focus. One of the champions of this approach, David Donoho, coined the term ''{\it compressed sensing\/}`` to emphasize the fact that $\ell_1$-minimization is not just a new way of massaging a ''complete`` set of measurements into a compact form, but rather a new way of thinking about how to measure things in the first place~\cite{Don06}.

This new way of thinking has profoundly practical implications. Making measurements can be expensive, in terms of time, money, or (in the case of, say, x-rays) damage done to the object being imaged. Compressive sensing has the potential to provide substantial cost savings without sacrificing accuracy. In one impressive numerical experiment, Cand\`{e}s, Romberg, and Tao~\cite{CanRT06} showed that a $512 \times 512$-pixel test image, known as the Logan-Shepp phantom, can be reconstructed exactly from 512 Fourier coefficients sampled along $22$ radial lines--with, in other words, more than $95\%$ of the ostensibly relevant data missing (see Figure 1).

A host of practical applications are now being explored, including new sensing techniques, new 
analog-to-digital converters, and a new digital camera with a single photon 
detector, being developed by Kevin Kelly, Richard Baraniuk and the Digital Signal Processing group at Rice (dsp.rice.edu/cs/cscamera)~\cite{TakLWDBSKB06,WakLDBSTKB06}.

\section{Mathematical foundations}

\subsection{Sparsity and undersampling}

The celebrated Nyquist-Shannon-Whittaker sampling theorem shows that a signal with bandwidth $2 \Omega$ is completely determined
by its uniform samples if and only if the samples are taken at least at the \emph{Nyquist rate} $\Omega/\pi$. This principle used
to underlie all signal acquisition techniques used in practice, such as consumer electronics, medical imaging, analog-to-digital
conversion and so on.  Compressive sampling puts forward a novel sampling paradigm that replaces the notion of band-limited signals 
with that of sparse signals. This new notion allows for dramatically ``undersampled'' signals to be captured and manipulated
using a very small amount of data.  The point of this section is to explain the basic mathematics behind this new theory.

Suppose $x$ is an unknown vector in $\R^N$ (a digital image or signal). We plan to sample $x$ using $n$ linear functionals of $x$ 
and then reconstruct. We are interested in the case  $n\ll N$, when we have many fewer measurements than the dimension of the signal
space. Such situations arise in many applications. For example, in biomedical imaging, far fewer measurements are typically collected
than the number of pixels in the image of interest. Further examples are provided by virtually any domain of science or technology 
where amounts of data are very large and costs of observation/acquisition/measurement are nontrivial.

The measurements $y_k$ are obtained by sensing $x$ against $n$ vectors $\phi_k\in \R^N$. Thus
$ y_k = \langle x, \phi_k \rangle$ for $k=1, \ldots, n$, or, equivalently 
\begin{equation} y = \Phi x \label{main_sys} \end{equation}
for some $n\times N$ \emph{measurement/sensing matrix} $\Phi$.  Thus we arrive at an underdetermined system
of linear equations, which, as is well known, in general has infinitely many solutions, so our problem 
is ill-posed. But suppose that our signal $x$ is \emph{sparse} or \emph{compressible}, i.e., that
is (essentially) depends only on a small number of degrees of freedom. To give a first impression of
the theory, we in fact assume that the signal can be written exactly as a linear combination of only a few basis 
vectors.

Mathematically the problem can be formulated as follows. Given a matrix $\Phi\in \R^{n\times N}$ with many 
more columns than rows ($n\ll N$), and a vector $y\in \R^n$, find a vector $x\in \R^N$ with a
minimum possible number of nonzero entries, i.e.,
\begin{equation}\label{ell0Prob1}
\text{minimize } \norm{x}{0} \;\; \text{subject to } \; \Phi x=y 
\end{equation}
where $\norm{x}{0}$ is the number of nonzero  entries of $x$ \cite{CorDIM03}. By allowing noise 
($\varepsilon \geq 0$), we obtain a variation of the problem (\ref{ell0Prob1}):
\begin{equation}\label{ell0Prob2}
\text{minimize }  \norm{x}{0} \;\;
\text{subject to } \; \norm{\Phi x-y}{2} \leq \varepsilon,
\end{equation}
These problems \emph{per se} are NP-hard even for $\varepsilon=0$,
see~\cite{GarJ79,Nat95}.

The classical, well studied, approach would be to minimize the $2$-norm $\norm{x}{2}$ in the
above problems, but this usually yields a solution vector $x$ that is full,
while for a sparse representation we would like to find a vector~$x$ with
few nonzero entries. 

The main approach taken in compressive sensing is to minimize the $1$-norm $\norm{x}{1}$ instead. 
\begin{equation}\label{ell1Prob1}
\text{minimize }  \norm{x}{1} \;\; 
\text{subject to }  \Phi x=y,
\end{equation}
and
\begin{equation}\label{ell1Prob2}
\text{minimize }   \norm{x}{1} \;\;
\text{subject to }  \norm{\Phi x-y}{2} \leq \varepsilon,
\end{equation}
respectively~\cite{CheDS99}, where $\norm{x}{1}\eqbd \sum_{i} \abs{x_i}$ (See Figure~2).
\begin{figure}\label{Lpnormfig1}
\begin{center}
\includegraphics[angle=0, width=0.5\textwidth]{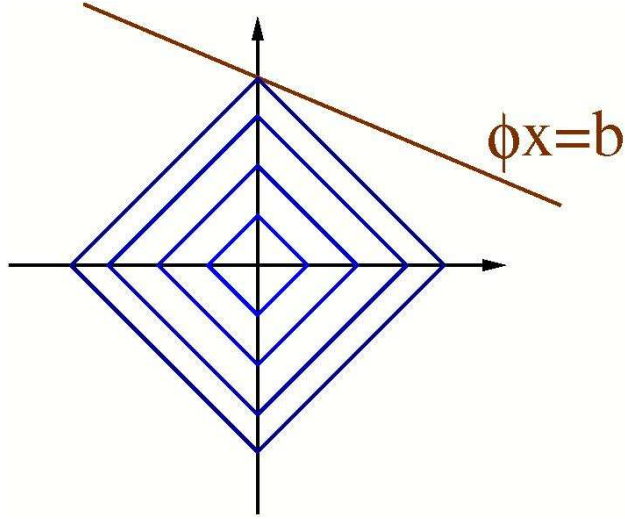}
\end{center}
\caption{$\ell_1$ minimization}
\end{figure}
Very surprisingly,  the $\ell_1$ minimization yields the same result as the $\ell_0$ minimization
in many cases of practical interest. This phenomenon was initially observed by engineers and
geophysicists, most notably  Claerbout and  Logan as early as  1970s (see \cite{TayBM79}), and  by Santosa 
and Symes in 1986 \cite{SanS86}, as mentioned in the introduction. In the last five years or so, a series of papers
\cite{CanRT06,Don0,Don1,DonH99,DonE03,Elad1,Fuchs1,GribNiel,ShaWM07} explained  why $\ell_1$ minimization
can recover sparse signals in a variety of practical setups.  In our next section, we give a few sample 
theorems about this remarkable phenomenon.

Finally, the $\ell_1$ minimization problem can be efficiently solved by convex programming or by linear 
programming (LP)~\cite{BoyV04}.  Most compressive sensing results due to Cand\`es, Donoho, Romberg, Tao and
others~\cite{Don0,Don1,CanRT06,Don06,CohDD06,CanT06} are all based on this method (see 
also~\cite{Zha05a,Zha05b,Zha05c}).   Other approaches include greedy algorithms, for instance, 
the so-called \emph{matching pursuit} introduced by Mallat and Zhang~\cite{MalZ93,PatRK93,Tro03,Tro04}. 
Recently many variations on matching pursuit have been proposed, among which are  orthogonal matching 
pursuit \cite{PatRK93,NeeV07}, stagewise orthogonal matching pursuit \cite{DonTIS06}, gradient pursuit 
\cite{BluD07}, and others.





\subsection{Incoherence and restricted isometry}

Given an $n\times N$ matrix $\Phi$, the first basic question is to determine 
whether $\Phi$ is good for compressive sensing, i.e., will lead to good 
recovery of sparse solutions to the equations $\Phi x =y$.  

Cand\`es and Tao \cite{Can06a}--\cite{CanT06b} introduced a necessary condition that guarantees
an estimate of its performance on classes of sparse vectors.  

\begin{definition}[{\cite{CanT05,CanRT06,CanR06}}] A matrix $\Phi$ is said to satisfy 
the {\em Restricted Isometry Property (RIP) of order $k$ with constant $\delta \eqbd 
\delta_k \in (0,1)$\/} if 
\begin{equation}
 (1-\delta_k) \|x \|_2^2 \leq \|  \Phi x \|_2^2  \leq (1+\delta_k) \|x \|_2^2     
\label{RIP}
\end{equation}
for any $x$ such that $\|x \|_0 \leq k$.
\end{definition}
 
It is straightforward to see that this condition can be reformulated as follows:
Consider $n\times \#T$ matrices $\Phi_T$ formed by the columns of $\Phi$ with indices
in the set $T$. Then the Gramian matrices 
$$ G_T\eqbd \Phi_T^t \Phi_T $$
are bounded and boundedly invertible on $l_2$ with bounds as in (\ref{RIP}), uniform
for all $T$ of size $\#T=k$. Since each matrix $G_T$ is symmetric and nonnegative 
definite, this is equivalent to each of these matrices having their eigenvalues 
in the interval $[1-\delta_k,1+\delta_k]$.
The role played by  RIP becomes clear from the following result of Cand\`es and Tao.

\begin{theorem}[\cite{CanRT06b,CanP06}]
 If the $n\times N$ matrix satisfies RIP of order $3k$ for some $\delta\in (0,1)$,
then, for any vector $x\in \R^N$, the $\ell_1$ minimization problem~(\ref{ell1Prob1}) 
has a solution $x^*$ such that
\begin{equation}
 \|  x- x^*  \|_2 \leq C \cdot { \|x-x_k \|_1 \over \sqrt{k} },
\label{RIP2}
\end{equation}
where $x_k$ denotes the best $k$-sparse approximation to $x$ and $C$ denotes a constant.
\end{theorem}

The condition~(\ref{RIP2}) means that the $\Phi$'s with higher
values of $k$ for which RIP is satisfied perform better in compressive sensing.
For example, if an $n\times N$ matrix $\Phi$ has the restricted isometry property of order $k$, 
then its performance in $l_2$-norm on the unit ball of $l_1^N$ is of order $C/\sqrt{k}$,
and the optimal performance is achieved if $\Phi$ satisfies RIP of order $k=\Theta(n/\log(N/n))$
\cite{DeVore07}. This is indeed achieved via various probabilistic constructions 
\cite{Don0,Don1,Don3,Don4,DonE04}.

A primary matrix measure related to RIP is {\em mutual incoherence\/} \cite{CheDS99,TroGMS03,Tro03,Tro04}: 
 $$ {\mathcal M}(\Phi)\eqbd \max_{i\neq j} |(\Phi^t \Phi)_{i,j}|, $$
i.e., the maximum inner product of distinct columns of $\Phi$. Since the
columns of $\Phi$ are usually normalized to be of $2$-norm $1$, the mutual
incoherence of a matrix is between $0$ and $1$. 

This notion can be generalized \cite{Jok07b} as follows:
For a given normalized matrix $\Phi\in \R^{n\times N}$, its {\it $k$-mutual
incoherence\/} ${\mathcal M}_k(\Phi)$ is defined by
$$ {\mathcal M}_k(\Phi)\eqbd \max_{\#S\leq k}\max_{i\neq j} |(\Phi_S^T \Phi_S)_{i,j}|. $$
The mutual incoherences ${\mathcal M}_k$ are intimately related to the best constant
$\delta_k$ with which the matrix $\Phi$ satisfies RIP of order $k$, but a full 
understanding of this connection has not been reached \cite{Jok07b,JokPfe06}.

A challenging aspect of RIP is its computational cost. Indeed, RIP is a property
of the submatrices of a specific size. At present, no subexponential-time algorithm 
is known for testing RIP.  Introducing other matrix measures may potentially help
in effectively verifying the RIP or finding other,  less demanding conditions for 
sparse recovery. 

One such weaker condition has been introduced by Cohen, Dahmen and DeVore in \cite{CohDD06}.
To motivate their condition, we first recall that a pair $(\Phi,\Delta)$ where 
$\Phi$ is a sensing matrix and $\Delta$ is a decoder, is called 
{\em instance-optimal of order $k$\/} for a normed space $(V, \|\cdot \|_V)$ if there
exists an absolute constant $C$ such that
$$   \|x-\Delta(\Phi x) \|_V \leq C \|x-x_k\|_V. $$
The matrix $\Phi$ has the {\em null space property\/} in $V$ if
$$  \| x\|_V \leq c \| x- x_k\|_V \qquad \hbox{\rm for all} \; x \text{ such that } \Phi x=0. $$
The importance of the null space property can be seen from the following result:

\begin{theorem}[{\cite{CohDD06, CohDD07}}]
Given an $n\times N$ matrix $\Phi$, a norm $\|\cdot \|_V$ and a value $k$, the instance
 optimality in $V$ with constant $C_0$ is equivalent to the null space property of $\Phi$
of order $2k$ with the constant $C_0/2$ in the sufficiency part and the same constant $C_0$
in the necessity part. 
\end{theorem}

Note that the null space property is preserved under row operations on the matrix $\Phi$
since, as its name suggests, it is simply a property of its null space. This property is
therefore less rigid than the RIP and may allow for a more efficient verification.

\subsection{Compressible signals}

In practice, most signals may not be exactly sparse in a given basis but may concentrate near a
sparse set. In fact, the most commonly used models in signal processing assume that the coefficients
of the signal with respect to, say, a wavelet basis, decay rapidly away from their essential support.
Smooth signals, images with bounded variation and those with bounded Besov norm are known to be of that type. 

Given a nearly sparse signal $x$, denote by $x_k$ its best $k$-sparse approximation, i.e., 
the vector obtained by keeping the $k$ largest coefficients of $x$ and discarding the rest.
Cand\`es, Romberg and Tao \cite{CanRT06} showed that the initial signal can be recovered with 
error of order $ \| x - x_k \|_1/ \sqrt{k}$ whenever the sensing matrix satisfies RIP of 
order $4k$ and the RIP constants $\delta_{3k}$ and $\delta_{4k}$ are not too close to $1$.

\begin{theorem}[{\cite{CanRT06b}}] \label{thm_comp_recovery} 
Let $\Phi$ satisfy RIP of order $4k$ with $\delta_{3k}+ 3 \delta_{4k}<2$. Then, for any signal
$x$, the solution $x^*$ to~(\ref{ell1Prob1}) satisfies 
$$ \|  x^* - x \|_2  \leq C \cdot  { \|x - x_k \|_1 \over \sqrt{k} }, $$
with a well-behaved constant $C$. 
\end{theorem}

\noindent
A similar result holds \cite{CanRT06} for stable recovery from imperfect measurements, i.e., 
in the setting of problem~(\ref{ell1Prob2}). All together, this indicates that $\ell_1$ minimization
stably recovers the largest $k$ coeffients of a nearly $k$-sparse vector even in the presence of noise.
 
This result is in fact optimal for important classes of signals: Let $x$ belong to the weak-$\ell_p$ 
ball or radius $R$, i.e., let the decreasing rearrangement of its coefficients
$ |x|_{(1)} \geq |x|_{(2)} \geq \cdots \geq |x|_{(N)}$ satisfy the condition
$$ |x|_{(i)} \leq R \cdot i^{-1/p}, \qquad i=1,\ldots, N.$$
This can be shown to imply 
$$  \| x - x_k \|_2 \leq C \cdot R \cdot k^{1/2-1/p} \quad {\rm and} \;\; 
  \| x - x_k \|_1 \leq C \cdot R \cdot k^{1-1/p}   $$
for some constant $C$.  Moreover, for generic elements in weak-$\ell_p$, 
no better estimates are obtainable. In other words, $\ell_1$ recovery achieves
an approximation error roughly as small as the error obtained by deliberately selecting 
the $k$ largest coefficients of the signal.

\subsection{Good sensing matrices}

  Most sampling algorithms developed so far in compressive sensing are based on randomization 
\cite{CanR05,Don06}.
Typically, the sensing matrices are produced by taking i.i.d. random variables with some given probability
distribution and then normalizing their columns. Such matrices are  guaranteed to perform well with very high probability, i.e., with the failure rate
exponentially small in the size of the matrix  \cite{Don06}.   Following \cite{Can06a}, we mention three
random constructions that are by now standard.

\begin{description}

\item{\bf Random matrices with i.i.d. entries.} Consider the matrix $\Phi$ with entries
drawn independently at random from the Gaussian probability distribution with mean zero
and variance $1/n$. Then \cite{CanT06,Don1}, with overwhelming probability,  the $\ell_1$ 
minimization~(\ref{ell1Prob1}) recovers 
$k$-sparse solutions whenever
$$ k \leq \const\cdot  n / \log (N/n). $$

\item{\bf Fourier ensemble.} Let $\Phi$ be obtained by randomly selecting $n$ rows from the $N\times N$
discrete Fourier transform and renormalizing the columns so that they have $2$-norm $1$.
If the rows are selected at random, then \cite{CanT06} as above, with overwhelming probability, the 
$\ell_1$ minimization~(\ref{ell1Prob1}) recovers $k$-sparse vectors for 
$$ k \leq \const \cdot n / (\log N)^6. $$

\item{\bf General orthogonal ensembles.} Suppose $\Phi$ is obtained by selecting $n$ rows from
an $N\times N$ orthonormal matrix $U$ and renormalizing the columns to be of unit length.
In the rows are selected at random, then \cite{CanT06} $k$-sparse recovery by $\ell_1$ minimization
(\ref{ell1Prob1}) is guaranteed with overwhelming probability provided that
$$ k \leq \const \cdot {1 \over {\cal M}^2(U ) } {n \over (\log N)^6 }. $$
Note that the Fourier matrix $U$ satisfies ${\cal M}(U)=1$, so this is a generalization
of the Fourier ensemble.

\end{description}

The natural problem already being addressed by several authors is how to achieve robust 
deterministic constructions of good CS matrices. Tao in \cite{Tao07} points out 
the importance of this  problem, as well as its similarity to other 
derandomization problems from theoretical computer science and combinatorics.
Several deterministic constructions are currently known (see, e.g., \cite{DeVore07,Ind08}).
However, the performance of matrices provided by these deterministic constructions is not yet on a par
with that of matrices arising probabilistically. 

To give several examples, DeVore in \cite{DeVore07} 
proposes a construction of cyclic matrices using finite fields that satisfy RIP of order $k$ for 
$k\leq C\sqrt{n} \log n / \log (N/n)$, which falls short of the above-mentioned range $k\leq C n / \log (N/n)$
known for probabilistic constructions.
Indyk in \cite{Ind08} and  Xu and Hassibi in \cite{XuHas} propose another scheme for compressive sensing with 
deterministic performance guarantees based on bipartite expander graphs.  Another flavor of randomness is 
introduced  in \cite{BajHRWN07} where random Toeplitz matrices are constructed with entries drawn independently 
from a given probability distribution.

\subsection{Optimality and $n$-widths}

The performance of the best sensing matrices $\Phi$, which is presently achieved
by random matrices with probabilistic guarantees, yields recovery of $k$-sparse
vectors using $n$ samples (so that the matrix $\Phi$ is $n\times N$) provided that
$$ k \leq \const \cdot n /\log(N/n). $$
In particular, a $k$-sparse vector can be recovered, say, by random projections,
of dimension $O(k\cdot \log(N/k))$ \cite{Don0}.  

For signals $x$ in the weak-$\ell_p$ ball of radius $R$,  $\ell_1$ recovery gives 
the error \cite{Don06}
$$ \|x^* -x \|_2 \leq \const \cdot R \cdot (n/\log(N/n))^{-1/p+1/2}.$$
It turns out that this performances cannot be improved even by using
possibly adaptive sets of measurements and reconstruction algorithms.

The matter turns out to be closely related to the issue of the so-called Gelfand widths  
\cite{Pin85} known from approximation theory: For a class $\cal F$, let $E_n ({\cal F})$
be the best reconstruction error from $n$ linear measurements
$$ E_n ({\cal F}) \eqbd \inf \sup_{f \in {\cal F}} \| f - D(y) \|_2, \quad y=\Phi f, $$
 where the infimum is over all sets of $n$ linear functionals and all reconstruction
algorithms $D$. The error $E_n({\cal F})$ is essentially equal~\cite{Pin85} to
the Gelfand width of the class $\cal F$ defined as
$$ d_n({\cal F}) \eqbd \inf_V \{  \sup_{f\in {\cal F}} \| P_V f \| : \codim(V) < n \}, $$
where $P_V$ is the orthogonal projector on the subspace $V$. Gelfand widths are known
for many classes of interest. In particular, Kashin \cite{Kashin77}, Garnaev and Gluskin \cite{GarGlu84}
showed that the Gelfand widths for the weak-$\ell_p$ ball of radius $R$ satisfy
$$  c \cdot R \cdot \left( \log(N/n) +1 \over n  \right)^{-1/p+1/2} \leq d_n({\cal F})
 \leq C \cdot R \cdot  \left(  \log(N/n) +1 \over n  \right)^{-1/p+1/2}.  $$
for some universal constants $c$ and $C$. 

This shows that the recovery provided by compressive sensing techniques is in fact optimal
for weak-$\ell_p$ norms in spite of being completely non-adaptive \cite{CohDD06,Can06a}. 
This is one more indication of the great potential of compressive sensing in applications.

\section{Connections with other fields}

\subsection{Statistical estimation}

Cand\`es \cite{Can06a} and Donoho \cite{Don06} point out a number of connections
of compressive sensing with ideas from statistics and coding theory.  We briefly
mention main ideas here.

In statistical estimation, the signal is assumed to be measured with stochastic errors
$$ y = \Phi x + z$$
where $z$ is a vector of i.i.d. (independent identically distributed) 
random variables with mean zero and variance $\sigma^2$. Very often, $z$ is 
assumed to be Gaussian.  The problem is again to recover $x$ from $y$.

One seeks to design an estimator whose accuracy depends on the information
content of the object $x$. The \emph{Dantzig selector} \cite{CanT06b} estimates $x$ by 
solving the convex program 
$$ \text{minimize } \|\tilde{x} \|_1 \;\; \text{ subject to } \;\;
\sup_{i} |(\Phi^T r)_i | \leq \lambda \sigma    $$
for some $\lambda>0$, where $r$ is the \emph{residual} 
$r\eqbd y - \Phi \widetilde{x}$. 
These ideas are very close to the so-called \emph{lasso} approach \cite{Tib96,EfrHJT04}.

Analogously to $\ell_1$ minimization in compressive sensing, 
the Danzig selector was shown \cite{CanT06b} to recover sparse and compressible 
signals with the number of measurements much smaller than the dimension of $x$
and within a logarithmic factor of the ideal mean squared error one would only
achieve with an \emph{oracle} supplying perfect information which coordinates 
are nonzero and which are above the noise level. 

\subsection{Error-correcting codes}

In coding theory \cite{Fel02,Fel03a,Fel03b,Fel03c}, a vector $x$ is transmitted to a remote receiver. 
The information $x$ is encoded using an $n\times N$ matrix $C$ with $n\ll N$.
Gross errors may occur during transmission, so that a fraction of the entries
of $Cx$ is completely corrupted. The location and the damage done to those entries
are unknown. It turns out that a constant fraction of errors with arbitrary
magnitudes can still be corrected \cite{CanT05} by solving a suitable linear minimization
problem. In fact, known methods recover the vector $x$ exactly provided 
the fraction of the corrupted entries is not too big \cite{CanT06b,Can06a}.

\subsection{Frame theory}

The theory of compressive sensing matrices closely resembles the basic theory of frames \cite{Chr03,Dau92,Mal98,You80}. 
A countable collection of elements $\{ f_i\}_{i\in I}$ is a \emph{frame} for a  Hilbert
space $H$ if there exist constants $0<A\leq B <\infty$ (the {\em lower\/} and {\em
upper frame bound\/}) such that, for all $g\in H$, 
$$ A \|g \|_H^2 \leq \sum_{i\in I} |\langle g, f_i  \rangle|^2 \leq B \|g \|_H^2.     $$
A frame is called {\em tight\/} if the upper and lower bounds are the same $A=B$. 
A frame is {\em bounded\/} if $\inf_{i\in I} \|f_i \|_H > 0$ (the condition $\sup_{i\in I}
\|f_i \|_H <\infty$ follows automatically from the definition of a frame). A frame is 
{\em unit norm\/} if $\|f_i\|_H=1$ for all $i\in I$. If $\{f_i \}_{i\in I}$ is a frame only 
for its closed linear span, it is called a {\em frame sequence. \/} A family 
$\{ f_i\}_{i\in I}$ is a {\em Riesz basic sequence\/} for $H$ if it is a {\em Riesz basis\/}
for its closed linear span, i.e., if, for some constants $0<A\leq B <\infty$ and for all
sequences of scalars $\{c_i \}_{i\in I}$,
$$ A \sum_{i\in I} |c_i|^2   \leq \| \sum_{i\in I} c_i f_i  \|_H^2 \leq  B \sum_{i\in I} |c_i|^2. $$ 
The analogy with the restricted isometry property is obvious, however,  the latter
is imposed only on submatrices formed from the original matrix. 

This analogy must be worth pursuing in both directions, i.e., looking for applications of the theory and methodology
of compressive sensing to frames and vice versa. Randomization techniques from compressive sensing could be
of particular interest in attacking problems from frame theory  (cf. \cite{BDDW,Tropp}).

\section{Practical implications}

Compressive sensing, and more generally the possibility of efficiently capturing 
sparse and compressible signals using a relatively small number of measurements,
paves the way for a number of possible applications.

\begin{description}
\item{\bf Data acquisition.}   New physical sampling devices may be designed that directly
record discrete low-rate incoherent measurements of the analog signal. This should be
especially useful in situations where large collections of samples may be costly, difficult 
or impossible to obtain.

\item{\bf Data compression.} The sparse basis in which the signal is to be represented may be
unknown or unavailable. However, a randomly designed $\Phi$ is suitable for almost 
all signals. We stress that these protocols are nonadaptive to the signal and
simply require to correlate it with a small number of other fixed vectors.   

\item{\bf Inverse problems.} The measurement system may have to satisfy rigid constraints
such as in MR angiography and other MR setups, where $\Phi$ records a subset of the 
Fourier transform. However, if a sparse basis exists that is also incoherent with $\Phi$,
then efficient sensing is possible.  

\end{description}

\begin{figure}\label{Camerafig}
\begin{center}
\includegraphics[angle=0, width=0.65\textwidth]{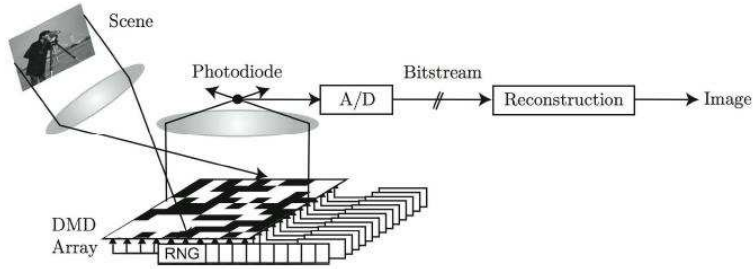}
\end{center}
\caption{The scheme of the CS camera}
\end{figure}

\noindent
A particularly interesting example of successful implementation of compressive sensing methodology
is provided by a digital camera newly developed by Richard Baraniuk and Kevin Kelly at Rice 
University (see dsp.rice.edu/cs/cscamera)~\cite{TakLWDBSKB06,WakLDBSTKB06}. 

In the detector array of a conventional digital camera, each pixel performs an analog-to-digital 
conversion; for example, the detector on a 5-megapixel camera produces 5 million bits for each image. 
This large amount of data is then dramatically reduced through a compression algorithms 
(using wavelet or other techniques) so as not to overburden typical storage and transfer capacities.

Rather than collect 5 million pixels for an image, the new camera  samples only a factor of 
about four times the 50,000 pixels that the jpg compression might typically output. These 
200,000 single-pixel measurements provide an immediate 25-fold savings in data collected 
compared with 5 megapixels.

The camera developed at Rice replaces the CCD array with a digital-micromirror device (DMD). 
A sequence of random projections is performed on the micromirror array, so that the image
``bounces off'' of each random pattern in the sequence, and the reflected light from each 
pattern is collected sequentially with a photodiode sensor that acts as the single-pixel 
detector (see Figure~3). After taking a sequence of essentially time-multiplexed measurements, 
a specific $\ell_1$ minimization algorithms decodes the picture out of the collected sequence 
of single-pixel measurements.

\end{document}